\documentclass{article}

\usepackage{natbib}

\bibpunct[:]{(}{)}{,}{a}{}{,}

\usepackage{times}

\usepackage{graphics}
\usepackage{graphicx}

\usepackage{amsmath,amssymb}
\usepackage{amsthm}
\usepackage{bm}
\usepackage{bbm}
\usepackage[normalem]{ulem}
\usepackage{relsize}

\usepackage[colorlinks=true,citecolor=blue,linkcolor=blue,breaklinks=true]{hyperref}

\usepackage{authblk}

\usepackage[hang]{footmisc}

\theoremstyle{definition}

\newtheorem{theorem}{Theorem}

\newtheorem{lemma}{Lemma}
\newtheorem{proposition}{Proposition}
\newtheorem{remark}{Remark}
\newtheorem{corollary}{Corollary}
\newtheorem{definition}{Definition}
\newtheorem{example}{Example}

\allowdisplaybreaks

\newcommand{\h}[1]{\mathbb{H}_T(#1)}
\newcommand{\hd}[1]{\mathbb{H}_T^{\dagger}(#1)}
\newcommand{\z}[1]{\mathbb{Z}_T(#1)}
\newcommand{\zd}[1]{\mathbb{Z}^{\dagger}_T( #1 )}
\newcommand{\zdt}[1]{\widetilde{\mathbb{Z}}^{\dagger}_T( #1 )}
\newcommand{\uu}{\mathbb{U}_T}
\newcommand{\uut}{\widetilde{\mathbb{U}}_T}
\newcommand{\jz}{\mathcal{J}^{(0)}}
\newcommand{\jo}{\mathcal{J}^{(1)}}
\newcommand{\thz}{\hat{\theta}_T^{(0)}}
\newcommand{\tho}{\hat{\theta}_T^{(1)}}
\newcommand{\tha}{\hat{\theta}_T}
\newcommand{\uhz}{\hat{u}^{(0)}_T}
\newcommand{\uho}{\hat{u}^{(1)}_T}
\newcommand{\uha}{\hat{u}_T}
\newcommand{\utz}{\tilde{u}^{(0)}_T}
\newcommand{\uto}{\tilde{u}^{(1)}_T}
\newcommand{\uta}{\tilde{u}_T}
\newcommand{\utinfz}{\tilde{u}_{\infty}^{(0)}}
\newcommand{\utinfo}{\tilde{u}_{\infty}^{(1)}}
\newcommand{\utinfa}{\tilde{u}_{\infty}}
\newcommand{\ta}{\tilde{a}_T}
\newcommand{\thom}{\underline{\hat{\theta}}_T^{(1)}}
\newcommand{\tsom}{\underline{\theta}^{*(1)}}

\begin{document}

\title{Penalized quasi likelihood estimation for variable selection
\footnote{
This work was in part supported by 
Japan Science and Technology Agency CREST JPMJCR14D7; 
Japan Society for the Promotion of Science Grants-in-Aid for Scientific Research 
No. 17H01702 (Scientific Research);  
and by a Cooperative Research Program of the Institute of Statistical Mathematics. 
}
}

\author[1,2,3]{Yoshiki Kinoshita}
\author[1,2,4]{Nakahiro Yoshida}
\affil[1]{Graduate School of Mathematical Sciences, University of Tokyo
\footnote{Graduate School of Mathematical Sciences, University of Tokyo: 3-8-1 Komaba, Meguro-ku, Tokyo 153-8914, Japan.}}
\affil[2]{CREST, Japan Science and Technology Agency}

\maketitle
\footnotetext[3]{E-mail: yoshikik@ms.u-tokyo.ac.jp}
\footnotetext[4]{E-mail: nakahiro@ms.u-tokyo.ac.jp}

\begin{abstract}
Penalized methods are applied to quasi likelihood analysis for stochastic differential equation models. In this paper, we treat the quasi likelihood function and the associated statistical random field for which a polynomial type large deviation inequality holds. Then penalty terms do not disturb a polynomial type large deviation inequality. This property ensures the convergence of moments of the associated estimator which plays an important role to evaluate the upper bound of the probability that model selection is incorrect. \vspace{6pt} \\
{\it keywards}: Quasi likelihood analysis; Polynomial type large deviation inequality; Variable selection.
\end{abstract}

\section{Introduction}
\label{intro}
Regularization methods, that impose a penalty term on a loss function, 
provide a tool for variable selection.
The method is useful because it performs estimation and variable selection simultaneously. 
Penalized estimators are generally expressed in the following form
\begin{align}
	\hat{\theta}_{\rm penalty} \in \underset{\theta\in\overline{\Theta}}{\rm argmin} \bigl\{ -L_n(\theta) + p(\theta) \bigr\}, \notag
\end{align}
where $\Theta$ is a parameter space, $L_n$ is a log likelihood function or $-L_n$ is equal to the sum of squared residuals and $p$ is a penalty term.
One of the most simple regularization methods is the Bridge \cite[]{frank1993statistical} that imposes the penalty term
\begin{align}
	p^{\rm Bridge}_{\lambda}(\theta) = \lambda\displaystyle\sum_{i = 1}^{\sf p}|\theta_i|^q \notag
\end{align}
on the least square loss function, where $q > 0$ is a constant, ${\sf p }$ is a dimension of an unknown parameter $\theta$ and $\lambda > 0$ is a tuning parameter.
For $q \leq 1$, the estimator performs variable selection.
Especially, when $q = 1$, the estimator is called the Lasso \cite[]{tibshirani1996regression}. Other than Bridge, various regularization methods have been proposed,
e.g. the smoothly clipped absolute deviation \cite[SCAD;][]{fan2001variable} and the minimax concave penalty \cite[MCP;][]{zhang2010nearly}. 
These methods are widely studied and extended in the regression analysis.
\cite{knight2000asymptotics} derived a $\sqrt n$-consistency of the Bridge estimator $\hat{\theta}_{\rm Bridge}$ and studied the limit distribution of $\sqrt{n}(\hat{\theta}_{\rm Bridge}-\theta^*)$ where $\theta^*$ is the true value of $\theta$.
\cite{zou2006adaptive} proposed the adaptive Lasso and derived its oracle property.
These results clarified the advantage of Bridge estimator with $q<1$ and adaptive Lasso estimator compared to Lasso estimator $\hat{\theta}_{\rm Lasso}$ in the sense of asymptotic efficiency because the limit distribution of $\sqrt{n}(\hat{\theta}_{\rm Lasso}-\theta^*)$ has a redundant term.

Applications of regularization methods to the quasi likelihood analysis (QLA) for stochastic models have been recently studied.
The penalized quasi maximum likelihood estimator is defined by
\begin{align}
	\hat{\theta}_T \in \underset{\theta\in\overline{\Theta}}{\rm argmin}\bigl\{ -\h{\theta} + p(\theta)\bigr\} \label{difofpestimator}
\end{align}
for a given quasi likelihood function $\mathbb{H}_T$ in these situations.
These approaches works well for various kinds of quasi likelihood functions; see e.g. \cite{belomestny2018low},
\cite{de2012adaptive} and \cite{gaiffas2019sparse}.
\cite{masuda2017moment} studied the moment convergence of the Lasso estimator under more general settings.
They derived the polynomial type large deviation inequality (PLDI) for the $L^1$-penalized contrast functions under suitable conditions. PLDI is an inequality given by \cite{yoshida2011polynomial}, which evaluate the random field
\begin{align}
	\z{u} = \exp\{\h{\theta^*+a_Tu}-\h{\theta^*}\}, \notag
\end{align}
where $a_T \in {\rm GL}({\sf p})$ is a deterministic sequence in the general linear group over $\mathbb{R}$ of degree ${\sf p}$.
This inequality plays a crucial role in QLA because it implies the uniform boundedness
\begin{align}
	\sup_{T>0}E[|a_T^{-1}(\hat{\theta}_T-\theta^*)|^m] <\infty \notag
\end{align}
and moment convergence
\begin{align}
	E[|a_T^{-1}(\hat{\theta}_T-\theta^*)|^{m'}] \rightarrow E[|\hat{u}_{\infty}|^{m'}] \notag
\end{align}
for some large $m>0$, every $m' \in (0,m)$ and a random variable $\hat{u}_{\infty}$ such that $a_T^{-1}(\hat{\theta}_T-\theta^*)\overset{d}{\rightarrow}\hat{u}_{\infty}$.
These properties are useful to investigate an asymptotic behavior of statistics which depends on the moment of $a_T^{-1}(\hat{\theta}_T-\theta^*)$; see e.g.
\cite{chan2011uniform},
\cite{shimizu2017moment},
\cite{suzuki2018penalized} and
\cite{umezu2019aic}.

PLDI can be derived form tractable conditions under the locally asymptotically quadratic (LAQ) settings \cite[]{yoshida2011polynomial}.
Actually, PLDI is obtained with LAQ on many kinds of models, e.g.
\cite{clinet2017statistical},
\cite{masuda2013convergence},
\cite{ogihara2014quasi}
and
\cite{uchida2013quasi}.

In this paper, we consider the quasi likelihood function $\mathbb{H}_T$ with LAQ and PLDI, and we study the penalized maximum likelihood estimator defined in (\ref{difofpestimator}).
Our penalty term can deal with many kinds of penalties including the Lasso, the Bridge and the adaptive Lasso.
The objective in this paper consists of two parts. One is to derive a polynomial type large deviation inequality for the penalized quasi likelihood random field and another is to study asymptotic behavior of the penalized quasi maximum likelihood estimators.

The rest of the paper is organized as follows.
Section \ref{spqle} describes our basic settings.
Section \ref{spldi} provides the polynomial type large deviation inequality for the penalized quasi likelihood function.
Section \ref{scvs} and \ref{sab} contain the basic property of the proposed estimator.
Section \ref{srvs} and \ref{smconv} give more advanced results, derived from the polynomial type large deviation inequality, for our estimator.
We apply our results to a stochastic differential model in Section \ref{sapp} and report the results of simulations in Section \ref{ssim}.

\section{Penalized quasi likelihood estimator} \label{spqle}

Let $\Theta$ be a bounded open set in ${\mathbb R}^{\sf p}$. 
We denote by $\theta^* \in \Theta $ the true value of an unknown parameter $\theta\in\Theta$. 
Given a probability space $(\Omega,{\cal F},P)$, we consider a sequence of random fields ${\mathbb H}_T:\Omega\times\overline{\Theta}\to{\mathbb R}$, $T\in\mathbb{T}$, where $\mathbb{T}$ is a subset of $\mathbb{R}_{\geq0}$ with $\sup \mathbb{T} = \infty$ and $\overline{\Theta}$ is a closure of $\Theta$. We assume $\h{\theta}$ is continuous for all $\omega\in\Omega$, where $\h{\theta}$ denotes the mapping $\overline{\Theta}\ni\theta\rightarrow\h{\theta,\omega}$ for each $\omega\in\Omega$. We call $\h{\theta}$ a quasi likelihood function and define the quasi maximum likelihood estimator (QMLE) $\hat\theta^{\rm QMLE}_T$ by
\begin{align}
	\hat\theta^{\rm QMLE}_T \in \underset{\theta\in\overline{\Theta}}{\rm argmax}~ \h{\theta}. \notag
\end{align}
Here we use this expression in the sense that $\hat\theta^{\rm QMLE}_T:\Omega\rightarrow\overline{\Theta}$ is a measurable mapping satisfying
\begin{align}
	\h{\hat\theta^{\rm QMLE}_T} = \underset{\theta\in\overline{\Theta}}{\rm max}~ \h{\theta} \notag
\end{align}
for all $\omega\in\Omega$.

If $\theta^*$ has sparsity, we can construct an estimator which performs parameter estimation and variable selection simultaneously by adding a penalty term to the quasi likelihood function. Let us consider the penalized quasi likelihood function 
\begin{align}
	\hd{\theta} = \h{\theta} - p_T(\theta) \label{290121-5}
\end{align}
and the penalized estimator
\begin{align}
	\tha \in \underset{\theta \in {\overline{\Theta}}}{\rm argmax}~\hd{\theta}, \notag
\end{align}
where $p_T:\overline{\Theta}\rightarrow\mathbb{R}_{\geq0}$ is a penalty function for every $T \in \mathbb{T}$.
In this paper, we assume that $p_T$ has the following expression
\begin{align}
	p_T(\theta) = \displaystyle\sum_{j = 1}^{\sf p}\xi_T^jp(\theta_j),\notag
\end{align}
where $\xi_T^j$ are (possibly random) positive sequences and $p:\mathbb{R}\rightarrow\mathbb{R}_{\geq 0}$ is a function satisfying $p(0) = 0$.
Indeed, this function is defined on $\mathbb{R}^{\sf p}$, but we consider the restriction to $\overline{\Theta}.$

In the following sections, we will denote
$\{ j; \theta^*_j = 0 \}$ and $\{ j; \theta^*_j \neq 0 \}$ by $\jz$ and $\jo$, 
respectively. Furthermore,  for a vector $x\in\mathbb{R}^{\sf p}$ and  a matrix $A\in\mathbb{R}^{{\sf p}\times {\sf p}}$, the vector $(x_j)_{j \in \mathcal{J}^{(k)}}$ and the matrix $(A_{ij})_{i \in \mathcal{J}^{(k)}, j \in \mathcal{J}^{(l)}}$ will be denoted by $x^{(k)}$ and $A^{(kl)}$, respectively, and we will express $x$ as $(x^{(0)},x^{(1)})$. We write $s(A,x) = Ax$, $s_j(A,x) = (Ax)_j$ and $s^{(k)}(A,x) = (Ax)^{(k)}$.
For tensors $A=(A_{i_1,...,i_d})$ and $B=(B_{i_1,...,i_d})$,
we denote $A[B]=\sum_{i_1,...,i_d}A_{i_1,...,i_d}B_{i_1,...,i_d}$.
Moreover we write
$A[u_1,\ldots,u_d] = A[u_1\otimes\cdots\otimes u_d]=\sum_{i_1,\ldots,i_d}A_{i_1,\ldots,i_d}u_1^{i_1}\cdots u_d^{i_d}$ for vectors $u_1=(u_1^{i_1})_{i_1},\ldots,u_d=(u_d^{i_d})_{i_d}$.
We denote by $u^{\otimes r}=u\otimes\cdots\otimes u$ the $r$ times tensor product of u.


\section{Polynomial type large deviation inequality}\label{spldi}

We make use of the quasi likelihood analysis (QLA) of \cite{yoshida2011polynomial} to examine the moment convergence of estimators for $\theta$ and to derive a central limit theorem of it.
Let $a_T \in{\rm GL}({\sf p})$ be a deterministic sequence satisfying $||a_T||\rightarrow 0$ as $T \rightarrow\infty$ and $\uu = \{ u \in \mathbb{R}^{\sf p}; \theta^* + a_Tu \in \overline{\Theta} \}$. Here $||A||$ denotes the square root of the maximum eigenvalue of $A'A$ for $A\in\mathbb{R}^{{\sf p}\times {\sf p}}$ and $A'$ is the transpose of $A$. 
Based on QLA, we define the random fields $\mathbb{Z}_T$ and $\mathbb{Z}^{\dagger}_T$ on $\uu$ by
\begin{align}
	\z{u}=\exp\biggl( \h{\theta ^*+a_Tu}-\h{\theta^*} \biggr)  \notag
\end{align}
and
\begin{align}
	\zd{u} = \exp\biggl( \hd{\theta^*+a_Tu}-\hd{\theta^*}  \biggr).\notag
\end{align}

Let $L > 0$ and $V_T(r) = \{u\in\uu; r \leq u\}$ for $r > 0$. We assume a polynomial type large deviation inequality (PLDI) in \cite{yoshida2011polynomial} for $\mathbb{Z}_T$.
\begin{description}
\item[[A1\!\!\!]] There exist constants $C_L > 0$ and $\varepsilon_L \in (0,1)$ such that\\
	\begin{align}
		P\Biggl[ \displaystyle\sup_{u\in V_T(r)} \z{u} \geq \exp (-r^{2-\varepsilon_L}) \Biggr] \leq \frac{C_L}{r^L} \label{pldi1}
	\end{align}
	 for all $r > 0, T > 0$.
\end{description}
Here the supremum of the empty set should read $-\infty$ by convention.
In this paper we assume that $a_T$ is a diagonal matrix, and write
\begin{align}
	a_T =
	\begin{pmatrix}
	\alpha^1_T&&&0\\
	&\alpha^2_T&&\\
	&&\ddots&\\
	0&&&\alpha^{\sf p}_T
	\end{pmatrix}.\notag
\end{align}
Note that $\prod_{1\leq j \leq {\sf p}}\alpha^j_T\neq0$ for all $T\in\mathbb{T}$.
Let $c_0$ be a positive constant. In order to estimate $\mathbb{Z}^{\dagger}_T$, we consider the following three conditions for the penalty term.
\begin{description}
	\item[[A2\!\!\!]] $p$ is differentiable except the origin.
\end{description}
\begin{description}
	\item[[A3\!\!\!]] For some positive constant $\varepsilon$,
	\begin{align}
	\displaystyle\sup_{-\varepsilon<x<\varepsilon}p(x)<\infty. \notag
	\end{align}
\end{description}
\begin{description}
	\item[[A4\!\!\!]] For all $j\in\jo$,
	\begin{align}
	\displaystyle\sup_{T\in\mathbb{T}} |\alpha_T^j\xi_T^j| \leq c_0 \notag
	\end{align}
	almost surely.
\end{description}

\begin{remark}
Given $\xi_T$, we can construct $\xi'_T$ satisfying [A4] by taking $\xi'_T$ such that
$\xi'^{j}_T=\min (\xi_T^j,(\alpha^j)^{-1}_Tc_0)$.
\end{remark}

\begin{example}[LASSO]\label{ex_lasso}
Define $\xi_T^j$ by $\xi_T^j = |\alpha_T^j|^{-1}$ and $p$ by $p(x)=|x|$, then the penalty term $p_T(\theta)=\sum_{j=1}^{\sf p}|\alpha_T^j|^{-1}|\theta_j|$ satisfies [A2]-[A4].
\end{example}

In the above setting, we can derive the PLDI for $\mathbb{Z}^{\dagger}_T$.

\begin{theorem}\label{tpldi}
	Given $L>0$, assume Conditions [A1]-[A4]. Then there exist constants $C'_L > 0$ and $\varepsilon'_L \in (0,1)$ such that\\
	\begin{align}
		P\Biggl[ \displaystyle\sup_{u\in V_T(r)} \zd{u} \geq \exp (-r^{2-\varepsilon'_L}) \Biggr] \leq \frac{C'_L}{r^L} \label{pldi}
	\end{align}
	 for all $r > 0, T > 0$.
\end{theorem}
\proof
By [A1], there exist constants $C_L > 0$ and $\varepsilon_L \in (0,1)$ satisfying (\ref{pldi1}) for all $r > 0, T > 0$.
Let $\varepsilon'_L\in(\varepsilon_L,1)$. For every $T>0$ and $r>0$, we have
\begin{align}
	P&\Biggl[ \displaystyle\sup_{u\in V_T(r)} \zd{u} \geq \exp (-r^{2-\varepsilon'_L}) \Biggr] \notag \\
	 &\leq P\Biggl[ \displaystyle\sup_{u\in V_T(r)} \zd{u}\exp\biggl\{\displaystyle\sum_{j\in\jz}\xi_T^jp\bigl(s_j(a_T,u)\bigr)\biggr\} \geq \exp (-r^{2-\varepsilon'_L}) \Biggr] \notag \\
	&\leq \displaystyle\sum_{n=0}^{\infty}P\left[ \sup_{\substack{2^nr \leq |u|\leq2^{n+1}r \\ u \in V_T(r)}} \z{u}\exp(B_1) \geq \exp (-r^{2-\varepsilon'_L}) \right], \notag
\end{align}
where
\begin{align}
	B_1 = -\displaystyle\sum_{j\in\jo}\xi_T^j\biggl[p\Bigl(\theta^*_j+s_j(a_T,u)\Bigr)-p(\theta^*_j)\biggr].\notag
\end{align}
Conditions [A2] and [A3] imply
\begin{align}
	\displaystyle\sup_{x \in U\setminus\{0\}} \frac{p(\theta+x)-p(\theta)}{x} < \infty \notag
\end{align}
for every $\theta \in \mathbb{R}\setminus\{0\}$ and every compact set $U\subset\mathbb{R}$.
Moreover, by definition of $\uu$, we observe that $\sup_{T\in\mathbb{T}}\sup_{u\in\uu}|a_Tu|<\infty$.
Therefore, from [A4], we have
\begin{align}	
	|B_1|&\leq \displaystyle\sum_{j \in \jo}\xi_T^j|s_j(a_T,u)|\Biggl|\frac{p(\theta^*_j+s_j(a_T,u))-p(\theta^*_j)}{s_j(a_T,u)}\Biggr| \notag \\
	&\leq c_0 K|u| \notag
\end{align}
for some  $K > 0$ which does not depend on $T$ and $r$. Then we have
\begin{align}
	\displaystyle\sum_{n=0}^{\infty}P&\left[ \sup_{\substack{2^nr \leq |u|\leq2^{n+1}r \\ u \in V_T(r)}} \z{u}\exp(B_1) \geq \exp (-r^{2-\varepsilon'_L}) \right] \notag \\
	& \leq  \displaystyle\sum_{n=0}^{\infty}P\left[ \displaystyle \sup_{\substack{2^nr \leq |u| \leq 2^{n+1}r \\ u \in V_T(r)}} \z{u} \geq \exp \left(-r^{2-\varepsilon'_L} -2^{n+1} c_0 Kr\right) \right] . \notag
\end{align}
Since $1 < 2-\varepsilon_L' < 2-\varepsilon_L$, there exists a constant $R_1 > 0$ such that
\begin{align}
	-r^{2-\varepsilon'_L} - 2^{n+1}c_0 Kr \geq -(2^nr)^{2-\varepsilon_L} \notag
\end{align}
for all $n \in \mathbb{N}$ and $r \geq R_1$.
By this inequality and [A1], we obtain
\begin{align}
	 \displaystyle\sum_{n=0}^{\infty}P&\left[ \displaystyle \sup_{\substack{2^nr \leq |u| \leq 2^{n+1}r \\ u \in V_T(r)}} \z{u} \geq \exp \left(-r^{2-\varepsilon'_L} -2^{n+1} c_0 Kr\right) \right] \notag \\
	&\leq \displaystyle\sum_{n=0}^{\infty}P\left[ \displaystyle \sup_{\substack{|u| \geq 2^nr \\ u \in V_T(r)}} \z{u} \geq \exp \Bigl( -(2^nr)^{2-\varepsilon_L} \Bigr) \right] \notag \\
	 &\leq \displaystyle\sum_{n=0}^{\infty}\frac{C_L}{(2^nr)^L} =  \frac{1}{r^L}\frac{2^LC_L}{2^L-1}\notag
\end{align}
for every $r \geq R_1$. Let $C'_L = \max\{ R_1^L, \frac{2^LC_L}{2^L-1} \} $, we complete the proof.
\qed

\vspace{16pt}

Let $\uha = a_T^{-1}(\hat{\theta}_T-\theta^*)$ then
\begin{align}
	\uha\in\underset{u\in\uu}{\rm argmax}~\zd{u}. \notag
\end{align}
PLDI derives the $L^m$-boundedness of $\uha$ (Proposition 1 of \cite{yoshida2011polynomial}).

\begin{proposition}\label{uhattight}
Let $L>m>0$. Suppose that there exists a constant $C_L$ such that
\begin{align}
P\Biggl[ \displaystyle\sup_{u\in V_T(r)} \zd{u} \geq 1 \Biggr] \leq \frac{C_L}{r^L} \notag
\end{align}
for all $T>0$ and $r>0$.
Then it holds that
\begin{align}
\sup_{T>0}E[|\uha|^m]<\infty. \label{utight}
\end{align}
In particular, $\uha = O_p(1)$ as $T\rightarrow\infty$ (i.e., for every $\epsilon >0$, there exist $\mathcal{T}\in\mathbb{T}$ and $M > 0$ such that $P(|\uha|>M)<\epsilon$ for all $T\geq \mathcal{T}$), under Conditions [A1]-[A4].
\end{proposition}


\section{Consistency of variable selection}\label{scvs}

In this section, we will derive the selection consistency of $\tha$.
Let $q \in (0,1]$, we consider the conditions for $p$.

\begin{description}
\item[[A5\!\!\!]] There exists $\lambda>0$ such that
\begin{align}
\displaystyle \lim _{x \rightarrow 0}\frac{p(x)}{{|x|}^q} = \lambda . \notag
\end{align}
\end{description}
\begin{description}
\item[[A6\!\!\!]] For every $j\in\jz$,
\begin{align}
	(\xi_T^j)^{-\frac{1}{q}}|\alpha_T^j|^{-1} \overset{p}{\rightarrow}0 \notag
\end{align}
as $T\rightarrow\infty$.
\end{description}
LASSO penalty in Example \ref{ex_lasso} derives PLDI, however, it does not satisfy [A6]. We give another example for [A6].

\begin{example}[Bridge type]\label{ex_bridge}
	Let $q < 1$ and $q' \in (q,1]$.
	Define $\xi_T^j$ by $\xi_T^j = |\alpha_T^j|^{-q'}$ and $p$ by $p(x)=|x|^q$, then the penalty term $p(\theta)=\sum_{j=1}^{\sf p} |\alpha_T^j|^{-q'}|\theta_j|^q$ satisfies [A2]-[A6].
\end{example}

Let $\ta$ be a diagonal matrix in $\mathbb{R}^{{\sf p}\times{\sf p}}$ satisfying $(\ta)_{jj} = (\xi^j_T)^{-\frac{1}{q}}$ for $j \in \jz$ and $(\ta)_{jj} = a^j_T$ for $j\in\jo$. Denote $a_T^{-1}\ta$ by $G_T$.
\begin{description}
\item[[A7\!\!\!]]
For every $M>0$,
\begin{align}
	\displaystyle\sup_{\substack{u,v\in\uu \\ |u|,|v|<M \\ u\neq v}}\frac{|\mathbb{H}_T(\theta^*+a_Tu)-\mathbb{H}_T(\theta^*+a_Tv)|}{|u-v|^q}||G_T^{(00)}||^q \overset{p}{\rightarrow}0 \notag
\end{align}
as $T\rightarrow\infty$.
\end{description}

\begin{remark}
Condition [A6] implies that
\begin{align}
	||G_T^{(00)}|| \overset{p}{\rightarrow}0 \label{gconv}
\end{align}
as $T\rightarrow\infty$.
We usually assume [A6] to ensure this convergence in this paper.
\end{remark}
\begin{remark}
Condition [A7] is a technical one, however, we can derive it easily from the differentiability of $\mathbb{H}_T$.
\end{remark}

\begin{description}
\item[[A7$'$\!\!\!]] For some $R>0$, the following conditions hold:
\begin{description}
\item[(i)] For every $T\in \mathbb{T}$, $\mathbb{H}_T$ is almost surely thrice differentiable with respect to $\theta$ on $B=B_R(\theta^*,\Theta)=\{ \theta\in\Theta ;  |\theta-\theta^*|<R  \}$,
\item[(ii)] $||a_T||\partial_{\theta}\mathbb{H}_T(\theta^*) = O_p(1)$,
\item[(iii)] $||a_T||^2\displaystyle\sup_{\theta\in B}|\partial_{\theta}^2\mathbb{H}_T(\theta)|=O_p(1)$,
\item[(iv)] $||a_T||^2\displaystyle\sup_{\theta\in B}|\partial_{\theta}^3\mathbb{H}_T(\theta)|=O_p(1)$.
\end{description}
\end{description}

\begin{proposition}\label{phconti}
	Assume [A6] and [A7$'$], then [A7] holds.
\end{proposition}

\proof
Take $R'<R$ satisfying $\overline{B_{R'}(\theta^*)}=\{ \theta\in \mathbb{R}^{\sf p}; |\theta-\theta^*|\leq R' \} \subset\Theta$.
For $M>0$, there exists a constant $\mathcal{T}_M\in\mathbb{T}$ such that $\theta^*+a_Tu\in B_{R'}(\theta^*)$ for every $T>\mathcal{T}_M$ and $u\in\mathbb{R}$ satisfying $|u| < M$.
Therefore, by Taylor's theorem
\begin{align}
	\bigl|\mathbb{H}_T(\theta^*&+a_Tu)-\mathbb{H}_T(\theta^*+a_Tv)\bigr| \notag\\
	\leq & \bigl|\partial_{\theta}\mathbb{H}_T(\theta^*)[a_Tu]-\partial_{\theta}\mathbb{H}_T(\theta^*)[a_Tv]\bigr| \notag\\
	& +\biggl|\int_0^1(1-s)\partial_{\theta}^2\mathbb{H}_T(\theta^*+sa_Tu)[(a_Tu)^{\otimes 2}]ds \notag\\
	&\qquad\qquad - \int_0^1(1-s)\partial_{\theta}^2\mathbb{H}_T(\theta^*+sa_Tv)[(a_Tv)^{\otimes 2}]ds\biggr|\notag\\
	\leq & A_1 + A_2 + A_3.\notag
\end{align}
for every $T>\mathcal{T}_M$ and every $u,v\in\mathbb{R}$ satifying $|u|,|v|<M$, where
\begin{align}
	A_1 &= ||a_T||\cdot\bigl|\partial_{\theta}\mathbb{H}_T(\theta^*)\bigr|\cdot|u-v|, \notag\\
	A_2 &= \biggl|\int_0^1(1-s)\partial_{\theta}^2\mathbb{H}_T(\theta^*+sa_Tu)[(a_Tu)^{\otimes 2}]ds \notag\\
	&\qquad\qquad - \int_0^1(1-s)\partial_{\theta}^2\mathbb{H}_T(\theta^*+sa_Tu)[(a_Tv)^{\otimes 2}]ds\biggr|\notag\\
	\intertext{and}
	A_3 &= \biggl|\int_0^1(1-s)\partial_{\theta}^2\mathbb{H}_T(\theta^*+sa_Tu)[(a_Tv)^{\otimes 2}]ds \notag\\
	&\qquad\qquad - \int_0^1(1-s)\partial_{\theta}^2\mathbb{H}_T(\theta^*+sa_Tv)[(a_Tv)^{\otimes 2}]ds\biggr|.\notag
\end{align}
However, [A7$'$](ii) and (\ref{gconv}) implies
\begin{align}
	\displaystyle\sup_{\substack{u,v\in\uu \\ |u|,|v|<M \\ u\neq v}}\frac{A_1}{|u-v|^q}||G_T^{(00)}||^q \overset{p}{\rightarrow}0 \label{prop2a1}
\end{align}
as $T\rightarrow\infty$.
If $\theta^*+a_Tu\in B_{R'}(\theta^*)$, then
\begin{align}
	A_2\leq ||a_T||^2\cdot|u+v|\cdot|u-v|\cdot\displaystyle\sup_{\theta\in B}|\partial_{\theta}^2\mathbb{H}_T(\theta)|. \notag
\end{align}
Therfore [A7$'$](iii) and (\ref{gconv}) implies
\begin{align}
	\displaystyle\sup_{\substack{u,v\in\uu \\ |u|,|v|<M \\ u\neq v}}\frac{A_2}{|u-v|^q}||G_T^{(00)}||^q \overset{p}{\rightarrow}0 \label{prop2a2}
\end{align}
as $T\rightarrow\infty$.
From Taylor's theorem, we have
\begin{align}
	\partial_{\theta}^2\mathbb{H}_T&(\theta^*+sa_Tu)=
	\partial_{\theta}^2\mathbb{H}_T(\theta^*+sa_Tv) \notag\\
	&+\int_0^1\partial_{\theta}^3\mathbb{H}_T\bigl(\theta^*+sa_Tu+s'(sa_Tu-sa_Tv)\bigr)ds'[sa_Tu-sa_Tv]\notag
\end{align}
for every  $s\in [0,1]$, $T>\mathcal{T}_M$ and every $u,v\in\mathbb{R}$ satifying $|u|,|v|<M$.
Since $\theta^*+sa_Tu+s'(sa_Tu-sa_Tv)\in B_{R'}(\theta^*)\subset B$, it follows that
\begin{align}
	A_3 \leq ||a_T||^3\cdot|v|^2\cdot|u-v|\cdot \sup_{\theta\in B}|\partial_{\theta}^3\mathbb{H}_T(\theta)| \notag
\end{align}
for every $T>\mathcal{T}_M$ and every $u,v\in\mathbb{R}$ satifying $|u|,|v|<M$.
Therefore [A7$'$](iv) and (\ref{gconv}) implies
\begin{align}
	\displaystyle\sup_{\substack{u,v\in\uu \\ |u|,|v|<M \\ u\neq v}}\frac{A_3}{|u-v|^q}||G_T^{(00)}||^q \overset{p}{\rightarrow}0 \label{prop2a3}
\end{align}
as $T\rightarrow\infty$.
From (\ref{prop2a1}), (\ref{prop2a2}) and (\ref{prop2a3}), we have the desired result.
\qed

The following theorem ensures that $\tha$ enjoys the consistency of variable selection.
\begin{theorem}\label{tcvs}
	Assume Conditions [A5] and [A7]. If $\uha=O_p(1)$, then
	\begin{align}
		P\Bigl( \thz = 0 \Bigr) \rightarrow 1 \label{t1main}
	\end{align}
	as $T\rightarrow\infty$.
\end{theorem}
\proof
Let $\mathcal{S}^1_{T,M} = \{ |\uha|<M , (0,\uho)\in\uu \}$.
For $M>0$ and $T \in \mathbb{T}$, define $\mathcal{S}^2_{T,M}$ by
\begin{align}
	\mathcal{S}^2_{T,M} = \Bigl\{ |\uhz| < M,\displaystyle\sum_{j\in\jz}\xi^j_Tp\bigl(s_j(a_T,\uha)\bigr) \geq \frac{\lambda}{2}|(G_T^{(00)})^{-1}\uhz|^q \Bigr\} \notag
\end{align}
and define $\mathtt{C}_{T,M}$ by
\begin{align}
	\mathtt{C}_{T,M} = \displaystyle\sup_{\substack{u,v\in\uu \\ |u|,|v|<M \\ u\neq v}}\frac{|\mathbb{H}_T(\theta^*+a_Tu)-\mathbb{H}_T(\theta^*+a_Tv)|}{|u-v|^q}. \notag
\end{align}
By definition,
\begin{align}
	\mathbb{Z}^\dagger_T(\uhz,\uho) - \mathbb{Z}^\dagger_T(0,\uho)&=\exp \Bigl( \mathbb{H}^\dagger\bigl(\theta^*+a_T\uha\bigr) - \mathbb{H}^\dagger\bigl(\theta^*+a_T(0,\uho)\bigr)\Bigr)\notag\\
	&= \exp \Bigl( \mathbb{H}\bigl(\theta^*+a_T\uha\bigr) - \mathbb{H}\bigl(\theta^*+a_T(0,\uho)\bigr)\Bigr)\notag\\
	&\qquad \times\exp\Bigl(-\sum_{j\in\jz}\xi^j_Tp\bigl(s_j(a_T,\uha)\bigr)\Bigr).\label{diffofh}
\end{align}
Therefore from [A1], we have
\begin{align}
	P(\thz \neq 0) &\leq P\Bigl(\mathbb{Z}^\dagger_T(\uhz,\uho) \geq \mathbb{Z}^\dagger_T(0,\uho), \uhz \neq 0 , (0,\uho)\in\uu \Bigr) \notag \\
	 &\qquad\qquad+ P((\mathcal{S}^1_{T,M})^c) \notag \\
	& \leq P\Bigl( \mathtt{C}_{T,M}|\uhz|^q \geq \frac{\lambda}{2}|(G_T^{(00)})^{-1}\uhz|^q, \uhz \neq 0 \Bigr) \notag\\
	&\qquad\qquad\qquad + P((\mathcal{S}^1_{T,M})^c) + P((\mathcal{S}^2_{T,M})^c). \notag
\end{align}
Therefore it suffices to estimate the following three probabilities:
\begin{align}
	P_1 &:= P\Bigl(\mathtt{C}_{T,M}||G_T^{(00)}||^q\geq\frac{\lambda}{2}\Bigr) ,\notag \\
	P_2 &:= P((\mathcal{S}^1_{T,M})^c) \notag
\intertext{and}
	P_3 &:= P((\mathcal{S}_{T,M}^2)^c). \notag
\end{align}
However, by [A7] we have $P_1\rightarrow0$ as $T\rightarrow\infty$.
Take $R>0$ satisfying $B_R(\theta^*) = \{\theta\in\mathbb{R}^{\sf p};|\theta-\theta^*|<R \}\subset\Theta$.
Since $\uha=O_p(1)$, $|a_T\uha|<R$ for large $T$ with large probability, therefore $P((0,\uho)\in\uu)\rightarrow1$ as $T\rightarrow\infty$.
Moreover, from [A5], for every $\epsilon>0$, there exist constants $M>0$ and $\mathcal{T}\in\mathbb{T}$ such that 
$P_2+P_3<\epsilon$ for every $T>\mathcal{T}$.
\qed


\section{Limit distribution}\label{sab}

In this section, we consider the central limit theorem of $\tha$.
Let $\uta = \ta^{-1}(\hat{\theta}_T-\theta^*)$ and $\uut (= \uut(\omega)) = \{ u\in\mathbb{R}^{\sf p}; \theta^*+\ta u\in\overline{\Theta} \}$.
Define the random field $\widetilde{\mathbb{Z}}^{\dagger}_T$ on $\uut$ by
\begin{align}
	\zdt{u} = \exp\biggl( \hd{\theta^*+\ta u}-\hd{\theta^*} \biggr), \notag
\end{align}
then
\begin{align}
	\uta\in\underset{u\in \uut}{\rm argmax}~\zdt{u} \notag
\end{align}
and
\begin{align}
\uha = G_T\uta . \notag
\end{align}
For convenience, we extend $\widetilde{\mathbb{Z}}_T$ to $\mathbb{R}^{\sf p}$ so that the extension has a compact support and $\displaystyle\sup_{u\in\mathbb{R}^{\sf p}\setminus\uut}\zdt{u}\leq\max_{u\in\partial\uut}\zdt{u}$.
In order to describe the limit distribution of $\uta$, we consider the following two conditions.

\begin{description}
\item[[A8\!\!\!]]
For all $M>0$,
\begin{align}
	\sup_{\substack{u\in\uu \\ |u|<M}}\Bigl| \mathbb{H}_T(\theta^*+a_Tu)-\mathbb{H}_T(\theta^*) \Bigr| = O_p(1) \notag
\end{align}
as $T\rightarrow\infty$.
\end{description}

\begin{description}
\item[[A9\!\!\!]]
For all $M>0$,
\begin{align}
	\displaystyle\sup_{\substack{u,v\in\uu \\ |u|,|v|<M \\ u\neq v}}\frac{|\mathbb{H}_T(\theta^*+a_Tu)-\mathbb{H}_T(\theta^*+a_Tv)|}{|u-v|^q} = O_p(1) \notag
\end{align}
as $T\rightarrow\infty$.
\end{description}

\begin{proposition}
	Assume Condition [A7$'$] is fulfilled, then [A8] and [A9] hold.
\end{proposition}

\proof
Similar to the proof of Proposition \ref{phconti}.
\qed

\begin{theorem}
Assume [A1], [A5], [A6], [A8] and [A9]. If $\uha =O_p(1)$, then
\begin{align}
	\uta = O_p(1) \notag
\end{align}
as $T\rightarrow\infty$.
\end{theorem}

\proof
By assumption and the definition of $\uta$,
\begin{align}
	\uto=\uho=O_p(1), \label{u1tight}
\end{align}
therefore it suffices to prove that $\utz = O_p(1)$.
For $T\in\mathbb{T}$, $R>0$, $M > 0$ and $\varepsilon_L$ as in [A1], define $\mathcal{S}^1_{T,R,M}$ by 
\begin{align}
	\mathcal{S}^1_{T,R,M} = \{ |\uha|<R,|\z{0,\uho}| > \exp(-M^{2-\varepsilon_L}) \}. \notag
\end{align}
Moreover define $P_1$, $P_2$ and $P_3(R)$ by
\begin{align}
	P_1 &= P\left(\displaystyle\sup_{\substack{|G_T^{(00)}u^{(0)}|\geq M \\ (G_T^{(00)}u^{(0)},\uho)\in\uu}} \zd{G_T^{(00)}u^{(0)},\uho} \geq \zd{0,\uho},\mathcal{S}^1_{T,R,M} \right),\notag \\
	P_2 &= P\left(\displaystyle\sup_{\substack{0<|G_T^{(00)}u^{(0)}| \leq M \\  (G_T^{(00)}u^{(0)},\uho)\in\uu}} \zd{G_T^{(00)}u^{(0)},\uho} \geq \zd{0,\uho},\mathcal{S}^1_{T,R,M} \right), \notag \\
\intertext{and}
	P_3(R) &= P\bigl((\mathcal{S}^1_{T,R,M})^c\bigr). \notag
\end{align}
Then for every $M > 0$ and $T\in\mathbb{T}$,
\begin{align}
	&\hspace{-5pt} P(|\utz|>M) \notag \\
	&\leq P\left(\displaystyle\sup_{\substack{|u^{(0)}|\geq M \\ (u^{(0)},\uto)\in\uut}} \zdt{u^{(0)},\uto} \geq \zdt{0,\uto} \right) \notag \\
	& \leq P_1 + P_2 + P_3(R). \notag
\end{align}
By [A1],
\begin{align}
	P_1 &\leq P\left(\displaystyle\sup_{\substack{|G_T^{(00)}u^{(0)}|\geq M \\ (G_T^{(00)}u^{(0)},\uho)\in\uu}} \notag \z{G_T^{(00)}u^{(0)},\uho} \geq \z{0,\uho},\mathcal{S}^1_{T,R,M} \right) \\
	&\leq P\Biggl[ \displaystyle\sup_{u\in V_T(M)} \z{u} \geq \exp (-M^{2-\varepsilon_L}) \Biggr] \notag\\
	& \leq \frac{C_L}{M^L} \label{t5p1}
\end{align}
for every $R>0$, $M>0$ and $T\in\mathbb{T}$.
Take $R'>0$ satisfying $\overline{B_{R'}(\theta^*)} \subset\Theta$ and take $\mathcal{T}_R$ satisfying $\sup_{T>\mathcal{T}_R}||a_T||R<R'$,
then for every $R>0$, $M>0$ and $T>\mathcal{T}_R$, $(0,\uho)\in\uu$ on $\mathcal{S}^1_{T,R,M}$.
Similarly to (\ref{diffofh}), by definition of $\mathtt{C}_{T,R}$, for every $u^{(0)}$ and $\uho$ such that $(0,\uho)$ and $(G_T^{(00)}u^{(0)},\uho)$ belong to $\uu$,
\begin{align}
	\zd{G_T^{(00)}u^{(0)},\uho}& - \zd{0,\uho} \leq \notag\\
	& \exp\Bigl(\mathtt{C}_{T,R}|G^{(00)}_Tu^{(0)}|^q
	 -\displaystyle\sum_{j\in\jz}\xi_T^jp\bigl((\xi_T^j)^{-\frac{1}{q}}u_j\bigr)\Bigr). \notag
\end{align}
Denote $B_1 = -\displaystyle\sum_{j\in\jz}\xi_T^jp\bigl((\xi_T^j)^{-\frac{1}{q}}u_j\bigr)$.
For $M>0$ and $T\in\mathbb{T}$, define $\mathcal{S}^2_{T,M}$ by
\begin{align}
	\mathcal{S}^2_{T,M} = \Biggl\{ \displaystyle\sup_{0<|G^{(00)}_Tu^{(0)}|\leq M}\frac{B_1}{|u^{(0)}|^q}<-\frac{\lambda}{2} \Biggr\}.\notag
\end{align}
Then we have
\begin{align}
	P_2 &\leq P\Bigl( \sup_{\substack{0<|G_T^{(00)}u^{(0)}| \leq M \\  \bigl(G_T^{(00)}u^{(0)},\uho\bigr)\in\uu}} (\mathtt{C}_{T,R}||G_T^{(00)}||^q-\frac{\lambda}{2}\bigr)|u^{(0)}|^q > 0 \Bigr) 
	+P\bigl((\mathcal{S}^2_{T,M})^c\bigr) \notag \\
	&= P\Bigl(  \mathtt{C}_{T,R}||G_T^{(00)}||^q > \frac{\lambda}{2}  \Bigr)
	+P\bigl((\mathcal{S}^2_{T,M})^c\bigr)\notag
\end{align}
for every $R>0$, $M>0$ and $T>\mathcal{T}_R$.
By [A9] and (\ref{gconv}), for every $\delta>0$, $R>0$ and $M>0$, there exists a constant $T_1(\delta,R,M) > \mathcal{T}_R$ such that
\begin{align}
	P\Bigl(  \mathtt{C}_{T,R}||G_T^{(00)}||^q > \frac{\lambda}{2}  \Bigr) < \delta \notag
\end{align}
for every $T\geq T_1(\delta,R,M)$.
Moreover [A5] implies that for every $\delta > 0$ and $M > 0$, there exists a constant $T_2(\delta,M)$ such that
\begin{align}
	P\bigl((\mathcal{S}^2_{T,M})^c\bigr) < \delta \notag
\end{align}
for every $T>T_2(\delta,M)$.
Therefore, for every $\delta>0$, $R>0$ and $M>0$, there exists a constant $T_3(\delta,R,M)$ such that
\begin{align}
	P_2  < 2 \delta \label{t5p2p3}
\end{align}
for every $T>T_3(\delta,R,M)$.
From (\ref{u1tight}), for every $\delta>0$, there exist $R_1>0$ and $T_4>0$ such that
\begin{align}
	P(|\uha|\geq R_1) < \delta \notag
\end{align}
for every $T>T_4$.
Moreover, [A8] implies that for every $\delta > 0$, there exist constants $T_5>0$ and $M_1>0$ such that
\begin{align}
	P_3(R_1) < \delta \label{t5p4}
\end{align}
for every $T>T_5$ and $M>M_1$.
We have the desired result from (\ref{t5p1}), (\ref{t5p2p3}) and (\ref{t5p4}).
\qed

\vspace{8pt}
We write $B(R) = \{u\in\mathbb{R}^{\sf p}; |u|\leq R\}$.
In order to describe the limit distribution of $\uta$, we introduce the local asymptotic quadraticity of $\mathbb{H}_T$.
\begin{definition}
The family $\mathbb{H}_T$ is called locally asymptotically quadratic (LAQ) at $\theta^*$ if there exist random vectors  $\Delta_T, \Delta \in \mathbb{R}^{\sf p}$,  random matrices $\Gamma_T, \Gamma \in \mathbb{R}^{{\sf p}\times{\sf p}}$ and random fields $r_T:\Omega\times\uu\rightarrow\mathbb{R}$ such that
\begin{description}
	\item[[A10\!\!\!]] (i) for every $T\in\mathbb{T}$ and $u\in\uu$
\begin{align}
	\mathbb{H}_T(\theta^*+a_Tu) - \mathbb{H}_T(\theta^*) = \Delta'_Tu - \frac{1}{2}u'\Gamma_T u + r_T(u), \notag
\end{align}
\begin{description}
	\item[(ii)] $\Gamma$ is almost surely positive definite,
	\item[(iii)] $(\Delta_T, \Gamma_T)\overset{d}{\rightarrow}(\Delta,\Gamma)$ as $T\rightarrow\infty$,
	\item[(iv)]  For all $R>0$, $\displaystyle\sup_{u\in B(R)}|r_T(u)|\overset{p}{\rightarrow}0$ as $T\rightarrow\infty$.
\end{description}
\end{description}
\end{definition}


\begin{remark}\rm 
One needs a certain global non-degeneracy of the random fields ${\mathbb H}_T$ as well as the LAQ property 
to prove the PLDI. 
Therefore [A1] is not redundant under [A10]. 
Moreover, the LAQ property will be used to identify the limit distribution of the estimators. 
\end{remark}

Let  
\begin{align}
{\mathbb Z}(u) =
\exp\bigg(\Delta'u-\frac{1}{2}u'\Gamma u\bigg)\qquad (u\in{\mathbb R}^{\sf p}) \notag
\end{align}
and let $\hat{C}({\mathbb R}^{\sf p}) = \{f\in C({\mathbb R}^{\sf p});\lim_{|u|\to\infty}|f(u)|=0\}$. Equip $\hat{C}({\mathbb R}^{\sf p})$ with the supremum norm. 
It is possible to extend ${\mathbb Z}_T$ from ${\mathbb U}_T$ to ${\mathbb R}^{\sf p}$ in such a way that 
the extended ${\mathbb Z}_T$ takes values in $\hat{C}({\mathbb R}^{\sf p})$ and 
$0\leq{\mathbb Z}_T(u)\leq \max_{v\in\partial{\mathbb U}_T} {\mathbb Z}_T(v)$ for all $u\in{\mathbb R}^{\sf p}\setminus{\mathbb U}_T$. 
We will write ${\mathbb Z}_T$ for the extended random field on ${\mathbb R}^{\sf p}$.

\begin{proposition}\label{290115-2}
Given $L>0$, suppose that [A1] and [A10] are fulfilled.
Let $m \in (0,L)$, then 
\begin{eqnarray*}
E\big[f\big(a_T^{-1}(\hat{\theta}_{\text{ML}}-\theta^*)\big)\big]
&\to& 
E[\big[f\big(\Gamma^{-1}\Delta\big)\big]
\end{eqnarray*}
as $T\to\infty$ for any continuous function $f:{\mathbb R}^{\sf p}\to{\mathbb R}$ satisfying $\displaystyle\limsup_{|u|\rightarrow\infty}|f(u)||u|^m<\infty$. 
\end{proposition}

\proof
The finite-dimensional convergence ${\mathbb Z}_T\to^{d_f}{\mathbb Z}$ is obvious. 
By [A10], we see that 
for any $\epsilon>0$, 
\begin{align}\label{290115-1}
\displaystyle\lim_{\delta\rightarrow\infty}\limsup_{T\rightarrow\infty}P\big(w_T(\delta,R)\geq\epsilon\big) = 0
\end{align}
where 
\begin{align}
w_T(\delta,R) = \displaystyle\sup_{\substack{u^1,u^2\in B(R)\\|u^1-u^2|\leq\delta}}\bigl|\log\mathbb{Z}_T(u^1)-\log\mathbb{Z}_T(u^2)\bigr|. \notag
\end{align}
Now the desired result follows from Theorem 4 of \cite{yoshida2011polynomial}. 
\qed

\begin{remark}\rm As a matter of fact, 
for Proposition \ref{290115-2}, 
the inequality of [A1] can be weakened. 
See \cite{yoshida2011polynomial} for details. 
\end{remark}


\begin{description}
\item[[A11\!\!\!]] For every $j\in\jo$, there exists a constant $\beta_j\in\mathbb{R}$ such that
\begin{align}
\displaystyle \xi^j_T\alpha^j_T\overset{p}{\rightarrow} \beta_j \notag
\end{align}
as $T\rightarrow\infty$.
\end{description}

\begin{example}
	Bridge type penalty  $p(\theta)=\sum_{j=1}^{\sf p} |\alpha_T^j|^{-q'}|\theta_j|^q$ as in Example \ref{ex_bridge} satisfies [A11].
	Especially, if $q' < 1$, then $\beta_j = 0$ for all $j\in\jo$.
\end{example}

Define the random field $\widetilde{Z}^\dagger$ on $\mathbb{R}^{\sf p}$ by
\begin{align}
	\widetilde{Z}^\dagger(u) = \exp\Bigl((\Delta^{(1)})'u^{(1)}& - \frac{1}{2}(u^{(1)})'\Gamma^{(11)}u^{(1)} \notag\\
	 &- \displaystyle \sum_{j\in\jz}\lambda|u_j|^q - \sum_{j\in\jo}\beta_j\frac{d}{dx}p(\theta^*_j)u_j\Bigr) \notag
\end{align}
then $\widetilde{Z}^\dagger$ has an unique maximizer $\utinfa = \displaystyle\underset{u\in\mathbb{R}^{\sf p}}{\rm argmax}~\widetilde{Z}^\dagger(u)$ where $\utinfz = 0$ and $\utinfo =$  \scalebox{0.92}{$(\Gamma^{(11)})^{-1}(\Delta^{(1)}-{\boldsymbol \psi}^{(1)})$}. Here ${\boldsymbol \psi}$ is some ${\sf p}$-dimensional vector such that \scalebox{0.92}{${\boldsymbol \psi}_j = \beta_j\frac{d}{dx}p(\theta^*_j)$} for $j\in\jo$.
In the above setting, we estimate the asymptotic distribution of $\uta$.
\begin{theorem}\label{asymptotic1}
Assume Conditions [A2], [A5], [A6], [A10] and [A11].
If $\uta = O_p(1)$, then
\begin{align}
	(\ta^{(0)})^{-1}(\thz-\theta^{*(0)}) \overset{p}{\rightarrow} 0 \notag
\end{align}
and
\begin{align}
	(\ta^{(1)})^{-1}(\tho-\theta^{*(1)}) \overset{d}{\rightarrow}  (\Gamma^{(11)})^{-1}(\Delta^{(1)}-{\boldsymbol \psi}^{(1)})  \notag
\end{align}
as $T \rightarrow \infty$.
\end{theorem}
\proof
It suffices to verify $\uta\overset{d}{\rightarrow}\utinfa$ as $T\rightarrow\infty$.
From [A2], [A5], [A6], [A10] and [A11], it follows that
\begin{align}
	\bigl(\zdt{u^1},\ldots,\zdt{u^n}\bigr)\overset{d}\rightarrow \bigl(\widetilde{Z}^\dagger(u^1),\ldots,\widetilde{Z}^\dagger(u^n)\bigr)\label{zconv}
\end{align}
as $T\rightarrow\infty$, for every $n\in\mathbb{N}$ and $u^1,\ldots,u^n\in \mathbb{R}$.\\
For $\delta>0$ and $R>0$, define $\tilde{w}_T(\delta,R)$ by
\begin{align}
	\tilde{w}_T(\delta,R)=\displaystyle\sup_{\substack{u^1,u^2\in B(R)\\|u^1-u^2|\leq\delta}}\bigl|\log\zdt{u^1}-\log\zdt{u^2}\bigr|. \notag
\end{align}
Then, from [A2], [A5], [A6], [A10] and [A11], we have
\begin{align}
	\displaystyle\lim_{\delta\rightarrow0}\limsup_{T\rightarrow\infty}P\bigl(\tilde{w}_T(\delta,R)>\epsilon\bigr)=0\label{modconti}
\end{align}
for each $R>0$ and $\epsilon>0$.
From (\ref{zconv}) and (\ref{modconti}), it follows that $\widetilde{\mathbb{Z}}^\dagger_T|_{B(R)}\overset{d}{\rightarrow}\widetilde{Z}^\dagger|_{B(R)}$ in $C(B(R))$ for every $R>0$, where $\widetilde{\mathbb{Z}}^\dagger_T|_{B(R)}$ and  $\widetilde{Z}|_{B(R)}$ denote the restriction of $\widetilde{\mathbb{Z}}^\dagger_T$ and $\widetilde{Z}^\dagger$ on $B(R)$ respectively.
Since $\uta = O_p(1)$, we have $\uta\overset{d}{\rightarrow}\utinfa$ as $T\rightarrow\infty$.
\qed

\begin{remark}
	The convergence of distribution of the result of Theorem \ref{asymptotic1} can be extended to the stable convergence, if we replace the convergence of distribution in [A10](iii) by the stable convergence.
\end{remark}

\section{Probability of variable selection}\label{srvs}

PLDI provides uniform boundedness of $\uha$ as mentioned in (\ref{utight}). It enable us to estimate a probability that a correct model is selected. Let $\eta \in (0,1]$. For $T\in\mathbb{T}$ and $R>0$, define $\mathsf{c}_{T,R}$ by
\begin{align}
	\mathsf{c}_{T,R} =  \displaystyle\sup_{\substack{u,v\in\uu \\ |a_Tu|,|a_Tv|<R||a_T||^{1-\eta} \\ u\neq v}}\frac{|\mathbb{H}_T(\theta^*+a_Tu)-\mathbb{H}_T(\theta^*+a_Tv)|}{|u-v|^q}. \notag
\end{align}
For $m>0$, we denote \scalebox{0.98}{$E[|\mathsf{c}_{T,R}|^{m}||G_T^{(00)}||^{qm}]$ by $c_T(m,R)$}.
If  \scalebox{0.98}{$E[|\mathsf{c}_{T,R}|^{m}||G_T^{(00)}||^{qm}]=$} $\infty$, we define $c_T(m,R)=\infty$.
\begin{remark}
The sequence $c_T(m,R)$ is expected to be small as $T\rightarrow\infty$. We will estimate it at the end of this section. 
\end{remark}
\begin{theorem}\label{trate}
	Given $m > 0$, suppose that the inequality (\ref{utight}) is fulfilled. Moreover assume [A5]. Then for every $R>0$ and $m_0>0$ there exists a positive constant $D_{m,m_0,R}$ such that
	\begin{align}
		P(\thz\neq 0) < D_{m,m_0,R}(||a_T||^{m\eta}+c_T(m_0,R)) \label{ratemain}
	\end{align}
	for every $T\in\mathbb{T}$.
\end{theorem}
\proof
By [A5], there exists a positive constant $R_1$ such that $p(x)>\lambda|x|^q/2$ for every $x$ satisfying $|x|<R_1$.
Take $R_2>0$ satisfying $\overline{B_{R_2}(\theta^*)}=\{ \theta\in \mathbb{R}^{\sf p}; |\theta-\theta^*|\leq R_2 \} \subset\Theta$ and let $R_3=\min \{ R_1, R_2 \}$.
For $T\in\mathbb{T}$ define $\mathcal{S}_T$ by
\begin{align}
	\mathcal{S}_T = \{ |a_T\uha| < R_3, |a_T\uha|<R||a_T||^{-\eta} \}. \notag
\end{align}
By (\ref{diffofh}) and definition of $R_3$ and $\mathsf{c}_{T,R}$,
\begin{align}
	P(\thz \neq 0) &\leq P\Bigl(\mathbb{Z}^\dagger_T(\uhz,\uho) \geq \mathbb{Z}^\dagger_T(0,\uho), \uhz \neq 0\Bigr) \notag \\
	& \leq P\Bigl(\mathsf{c}_T|\uhz|^q \geq \frac{\lambda}{2}|(G_T^{(00)})^{-1}\uhz|^q, \uhz \neq 0 \Bigr) + P(\mathcal{S}^c_T). \notag
\end{align}
Therefore it suffices to estimate the following two probabilities:
\begin{align}
	P_1 &:= P\Bigl(\mathsf{c}_T||G_T^{(00)}||^q\geq\frac{\lambda}{2}\Bigr) ,\notag
\intertext{and}
	P_2 &:= P(\mathcal{S}_T^c). \notag
\end{align}
By Markov's inequality,
\begin{align}
	P_1 \leq \Bigl(\frac{2}{\lambda}\Bigr)^{m_0} c_T(m_0,R) \label{t6p3}
\end{align}
and
\begin{align}
	P_2 \leq (R^{-m}||a_T||^{m\eta}+R_3^{-m}||a_T||^{m})\displaystyle\sup_{T\in\mathbb{T}}E[|\uha|^m] \label{t6p4}
\end{align}
for every $T\in\mathbb{T}$. 

From (\ref{utight}), (\ref{t6p3}) and (\ref{t6p4}), we have (\ref{ratemain}) for some $D_{m,m_0,R}$.

\qed

\vspace{8pt}
Theorem \ref{trate} gives an upper bound of the probability of overfitting, however, we need to estimate the probability of underfitting. Let $\thom=\min_{j\in\mathcal{J}^{(1)}}|\hat{\theta}_{T,j}|$.
\begin{theorem}\label{trate1}
	Given $m>0$, suppose that the inequality (\ref{utight}) is fulfilled. Then there exists a positive constant $D_m$ such that
	\begin{align}
		P(\thom = 0) < D_m||a_T||^m \label{rate1main}
	\end{align}
	for every $T\in\mathbb{T}$.
\end{theorem}
\proof
Let $\tsom=\min_{j\in\mathcal{J}^{(1)}}|\theta^*_j|$, then
\begin{align}
	P(\thom = 0)\leq P(|a_T\uha|\geq\tsom). \notag
\end{align}
Moreover, by Markov's inequality,
\begin{align}
	P(\thom = 0)\leq \frac{||a_T||^m}{|\tsom|^m}\displaystyle\sup_{T\in\mathbb{T}}E[|\uha|^m] \notag
\end{align}
for all $T\in\mathbb{T}$. Therefore from assuption, we have (\ref{rate1main}).
\qed 

\vspace{8pt}
We obtain the following corollary from Theorem \ref{trate} and Theorem \ref{trate1}:
\begin{corollary}\label{crate}
	Given $m>0$, suppose that the inequality (\ref{utight}) is fulfilled. Moreover assume that Condition [A5] holds. Then for every $R>0$ and $m_0>0$, there exists a positive constant $D_{m,m_0,R}$ such that
	\begin{align}
		P\bigl(\{j;\hat{\theta}_{T,j} = 0 \} \neq \jz\bigr) < D_{m,m_0,R}(||a_T||^{m\eta} + c_T(m_0,R)) \notag
	\end{align}
	for every $T\in\mathbb{T}$.	
\end{corollary}

\subsection{Estimation of $c_T(m_0,R)$}
In this subsection, we assume that $G_T$ is deterministic for simplicity and denote $G_T^{(00)}$ by $g_T$.
Let $R^*_0>0$ and $m_1>0$.
\begin{description}
\item[[A12\!\!\!]] {(i)}
For every $T\in \mathbb{T}$, $\mathbb{H}_T$ is almost surely thrice differentiable with respect to $\theta$ on $B^*=B_{R^*_0}(\theta^*,\Theta)=\{ \theta\in\Theta ;  |\theta-\theta^*|<R^*_0  \}$,
\begin{description}
\item[(ii)] $\sup_{T\in\mathbb{T}}E\Bigl[||a_T||^{m_1}|\partial_{\theta}\mathbb{H}_T(\theta^*)|^{m_1}\Bigr] < \infty$,
\item[(iii)] $\sup_{T\in\mathbb{T}}E\Bigl[||a_T||^{2m_1}\displaystyle\sup_{\theta\in B^*}|\partial_{\theta}^2\mathbb{H}_T(\theta)|^{m_1}\Bigr]<\infty$,
\item[(iv)] $\sup_{T\in\mathbb{T}}E\Bigl[||a_T||^{2m_1}\displaystyle\sup_{\theta\in B^*}|\partial_{\theta}^3\mathbb{H}_T(\theta)|^{m_1}\Bigr]<\infty$.
\end{description}
\end{description}

\begin{proposition}\label{prop_est_c_deterministic}
Assume that Condition [A12] holds.
Then there exist positive constants $R_0>0$ and $K>0$ such that
\begin{align}
	c_T(m_1,R_0)\leq K(||a_T||^{-\eta(2-q)}g_T^{q})^{m_1} \notag
\end{align}
for every $T\in\mathbb{T}$ satisfying $||a_T||\leq1$.
\end{proposition}

\proof
Take $R_0 \leq R^*_0$ satisfying that $\overline{B_{R_0}(\theta^*)}\subset\Theta$, then $\theta^*+a_Tu\in B^*$ for every $u$ satisfying that $|a_Tu|<R_0$.
Similarly to the proof of Proposition \ref{phconti}, we have
\begin{align}
	\bigl|\mathbb{H}_T(\theta^*&+a_Tu)-\mathbb{H}_T(\theta^*+a_Tv)\bigr| \leq A_1 + A_2 + A_3 \notag
\end{align}
for every $u, v\in\mathbb{R}$ satisfying $|a_Tu|, |a_Tv| < R_0$ where
\begin{align}
	A_1 &= |a_T(u-v)|\cdot|\partial_{\theta}\mathbb{H}_T(\theta^*)|,\notag\\
	A_2 &= |a_T(u+v)|\cdot|a_T(u-v)|\cdot\sup_{\theta\in B^*}|\partial^2_{\theta}\mathbb{H}_T(\theta)|\notag
	\intertext{and}
	A_3 &= |a_Tv|^2\cdot|a_T(u-v)|\cdot\sup_{\theta\in B^*}|\partial^3_{\theta}\mathbb{H}_T(\theta)|.\notag
\end{align}
If both $|a_Tu|$ and $|a_Tv|$ are less than $||a_T||^{1-\eta}R_0$, then
\begin{align}
	\frac{|a_T(u-v)|}{|u-v|^q} &\leq |a_T(u-v)|^{1-q}\frac{||a_T||^q|u-v|^q}{|u-v|^q} \notag\\
	&\leq (2R_0||a_T||^{1-\eta})^{1-q}||a_T||^q \notag\\
	&=(2R_0)^{1-q}||a_T||^{(1-\eta)(1-q)+q}. \notag
\end{align}
Therefore
\begin{align}
	E[|\mathsf{c}_{T,R_0}|^{m_1}] &\leq E\Biggl[\sup_{\substack{u,v\in\uu \\ |a_Tu|,|a_Tv|<R_0||a_T||^{1-\eta} \\ u\neq v}}
	\biggl(\frac{A_1+A_2+A_3}{|u-v|^q}\biggr)^{m_1}\Biggr] \notag\\
	 &\leq E\Biggl[\sup_{\substack{u,v\in\uu \\ |a_Tu|,|a_Tv|<R_0||a_T||^{1-\eta} \\ u\neq v}}
	\frac{3^{m_1-1}(A_1^{m_1}+A_2^{m_1}+A_3^{m_1})}{|u-v|^{qm_1}}\Biggr] \notag\\
	&\leq 3^{m_1-1}(A'_1+A'_2+A'_3)\notag
\end{align}
for every $T\in\mathbb{T}$ satisfying $||a_T||\leq1$ where
\begin{align}
	A'_1&=E\biggl[\Bigl((2R_0)^{1-q}||a_T||^{(1-\eta)(1-q)+q}|\partial_{\theta}\mathbb{H}_T(\theta^*)|\Bigr)^{m_1}\biggr],\notag\\
	A'_2&=E\biggl[\Bigl((2R_0)^{2-q}||a_T||^{(1-\eta)(2-q)+q}\sup_{\theta\in B^*}|\partial^2_{\theta}\mathbb{H}_T(\theta)|\Bigr)^{m_1}\biggr]\notag\\
	\intertext{and}
	A'_3&=E\biggl[\Bigl(2^{1-q}R_0^{3-q}||a_T||^{(1-\eta)(3-q)+q}\sup_{\theta\in B^*}|\partial^3_{\theta}\mathbb{H}_T(\theta)|\bigr)^{m_1}\biggr].\notag
\end{align}
If $||a_T||\leq1$, then Conditions [A12](ii), (iii) and (iv) imply
\begin{align}
	A'_1&\leq K'||a_T||^{-\eta(1-q)m_1}\leq K'||a_T||^{-\eta(2-q)m_1},\notag\\
	A'_2&\leq K'||a_T||^{-\eta(2-q)m_1}\notag\\
	\intertext{and}
	A'_3&\leq K'||a_T||^{(1-\eta(3-q))m_1}\leq K'||a_T||^{-\eta(2-q)m_1}\notag
\end{align}
for some $K'>0$, respectively.
Let $K > 3^{m_1}K'$, we have the desired result.
\qed


\begin{example}\label{exam_est_c}
Define $p(\theta)=\sum_{j=1}^{\sf p} |\alpha_T^j|^{-q'}|\theta_j|^q$ as in Example \ref{ex_bridge}, then we have  $g_T=||a_T||^{(q'-q)/q}$.
Let $m,m_1>0$ and suppose that the inequality (\ref{utight}) is fulfilled. Moreover assume that Conditions [A5] and [A12] hold.
Let $\eta=(q'-q)m_1/(m+2(1-q)m_1)$.
Then by Proposition \ref{prop_est_c_deterministic} and Corollary \ref{crate}, there exists a constant $D_{m,m_1}$ such that
\begin{align}
	P\bigl(\{j;\hat{\theta}_{T,j} = 0 \} \neq \jz\bigr) < D_{m,m_1}||a_T||^{\frac{(q'-q)mm_1}{m+2(1-q)m_1}} \notag
\end{align}
for every $T\in\mathbb{T}$.
\end{example}

\subsection{The case of random $G^{(00)}_T$}
Now we turn to the estimation of $c_T(m_0,R)$ in the case where $G_T^{(00)}$ is random.
Let $m_2>0$.
\begin{proposition}\label{prop_est_c_rand}
Assume that Condition [A12] holds. Then for every $m_2>0$, there exist positive constants $R_0>0$ and $K>0$ such that
\begin{align}
	c_T\Bigl(\frac{m_1m_2}{qm_1+m_2},R_0\Bigr) \leq K||a_T||^{-\frac{\eta(2-q)m_1m_2}{qm_1+m_2}}\bigl(E[||G_T^{(00)}||^{m_2}]\bigr)^{\frac{qm_1}{qm_1+m_2}} \notag
\end{align}
for every $T\in\mathbb{T}$ satisfying $||a_T||\leq1$.
\end{proposition}

\proof
Similarly to the proof of Proposition \ref{prop_est_c_deterministic}, there exist constants $R_0>0$ and $K'>0$ such that
\begin{align}
	E[|\mathsf{c}_{T,R_0}|^{m_1}] \leq K'||a_T||^{-\eta(2-q)m_1}\notag
\end{align}
for every $T\in\mathbb{T}$ satisfying that $||a_T|| \leq 1$.
Therefore from H\"{o}lder's inequality, we have
\begin{align}
	E[|\mathsf{c}_{T,R_0}|^{m_0}||G_T^{(00)}||^{qm_0}] &\leq \bigl(E[|\mathsf{c}_{T,R_0}|^{m_1}]\bigr)^{\frac{m_2}{qm_1+m_2}}\bigl(E[||G_T^{(00)}||^{m_2}]\bigr)^{\frac{qm_1}{qm_1+m_2}} \notag\\
	&\leq  K'^{\frac{m_2}{qm_1+m_2}}||a_T||^{-\frac{\eta(2-q)m_1m_2}{qm_1+m_2}}\bigl(E[||G_T^{(00)}||^{m_2}]\bigr)^{\frac{qm_1}{qm_1+m_2}} \notag
\end{align}
where $m_0=m_1m_2/(qm_1+m_2)$.
\qed

\begin{example}
In example \ref{exam_est_c}, we have $g_T=||a_T||^{(q'-q)/q}$.
Here we assume that there exist positive constants $m_2$ and $K'$ such that $E[||G_T^{(00)}||^{m_2}]\leq K'||a_T||^{(q'-q)m_2/q}$ for every $T\in\mathbb{T}$.
Suppose that Condition [A12] holds.
Then by proposition \ref{prop_est_c_rand}, there exists a constant $K>0$ such that
\begin{align}
	c_T \leq K||a_T||^{\frac{(q'+q\eta-2\eta-q)m_1m_2}{qm_1+m_2}} \notag
\end{align}
for every $T\in\mathbb{T}$ satisfying that $||a_T||\leq 1$.
Let $m>0$ and suppose that the inequality (\ref{utight}) is fulfilled.
Moreover assume that Condition [A5] holds.
Let $\eta=(q'-q)m_1m_2/(qmm_1+mm_2+2m_1m_2-qm_1m_2)$, then by Corollary \ref{crate}, there exists a constant $D_m$ such that
\begin{align}
	P\bigl(\{j;\hat{\theta}_{T,j} = 0 \} \neq \jz\bigr) < D_{m,m_1,m_2}||a_T||^{\frac{(q'-q)mm_1m_2}{qmm_1+mm_2+(2-q)m_1m_2}} \notag
\end{align}
for every $T\in\mathbb{T}$
\end{example}

\section{Moment convergence}\label{smconv}
In this section, we will study the moment convergence of $\uha$. The following theorem is a consequence of PLDI:
\begin{theorem}\label{tmconv1}
	Given $m>0$, suppose that (\ref{utight}) holds. Moreover assume that the conclusion of Theorem \ref{asymptotic1} holds. Then we have
\begin{align}
E[f(\uha)] \rightarrow E[f(\utinfa)] \notag
\end{align}
as $T\to\infty$ for any continuous function $f:{\mathbb R}^{\sf p}\to{\mathbb R}$ satisfying $\displaystyle\limsup_{|u|\rightarrow\infty}|f(u)||u|^{-m} = 0$. 
\end{theorem}
\proof
Condition (\ref{utight}) implies an uniform integrability of $\{f(\uha)\}_{T\in\mathbb{T}}$.
By assumption, $\uha\overset{d}{\rightarrow}\utinfa$ as $T\rightarrow\infty$.
Therefore we can obtain the desired result.
\qed

\vspace{8pt}
Theorem \ref{tmconv1} suggests that $\lim_{T\rightarrow\infty}E\bigl[|(a_T^{(00)})^{-1}\thz|^m\bigr] = 0$ for $m\in(0,L)$. From Theorem \ref{trate}, we have another estimation of $\thz$. Let $\Psi_T\in{\rm GL}(|\jz|)$ be a deterministic sequence of positive matrices. For $m>0$, we consider the condition
\begin{description}
\item[[A13\!\!\!]] $||\Psi_T||^{m^*}(||a_T||^{m\eta} + c_T(m_0,R))\rightarrow0$ as $T\rightarrow\infty$.
\end{description}
\begin{theorem}\label{tmconv2}
	Given $m,m^*,m_0>0$, suppose that Condition [A13] and the inequality (\ref{ratemain}) hold. Then
	\begin{align}
	E[|\Psi_T\thz|^{m^*}]\rightarrow0 \notag
	\end{align}
	as $T\rightarrow\infty$. \notag
\end{theorem}
\proof
Let $\Theta_{\max}$ = $\sup_{\theta\in\overline{\Theta}}|\theta|$, then
\begin{align}
	E[|\Psi_T\thz|^{m^*}] &\leq \Theta_{\max}^{m^*}||\Psi_T||^{m^*}P(\thz\neq0). \notag
\end{align}
Therefore, from [A13], we have the desired result.
\qed


\section{Application}\label{sapp}

In this section, we apply our results to stochastic differential models.
Consider a stochastic regression model specified by the 
stochastic integral equation
\begin{eqnarray} \label{sde1}
  Y_t &=& Y_0+\int_0^t b_s ds+ \int_0^t \sigma (X_s, \theta )dw_s, \quad t \in [0,T]. 
\end{eqnarray}
Here given a stochastic basis $(\Omega, {\cal F}, ({\cal F}_t)_{t \in [0,T]}, P)$, 
$w$ is an ${\sf r}$-dimensional standard Wiener process, and 
$b$ and $X$ are progressively measurable processes 
taking values in ${\mathbb R}^{{\sf m}}$ and ${{\mathbb R}}^{\sf d}$, respectively. 
The function 
$\sigma:{\mathbb R}^{\sf d} \times \Theta\to{\mathbb R}^{{\sf m}} \otimes {\mathbb R}^{\sf r}$ has 
an unknown parameter $\theta\in\Theta$, a bounded open set of ${\mathbb R}^{{\sf p}}$. 
Furthermore, if $b_t = b(Y_t,t)$ and $X_t =(Y_t,t)$, then $Y$ can be a time-inhomogeneous diffusion process. 
%
We want to estimate $\theta$ from the observations 
$(X_{t_j},Y_{t_j})_{j=0,...,n}$, $t_j = j h$ for $h=h_n = T/n$. 
No data of $b_t$ is available. 

High frequency data under finite time horizon will be treated, 
that is, $T$ is fixed and $n$ tends to $\infty$. 
This is a standard setting in finance. 
We will consider 
the penalized quasi likelihood analysis for the volatility parameter $\theta$. 
To apply the results in Sections \ref{spqle}-\ref{smconv}, 
we will use $n$ for $T$ of Section \ref{spqle},
while $T$ denotes the fixed terminal of the observations in what follows.



Let $S(x,\theta)=\sigma(x,\theta)^{\otimes2}=\sigma(x,\theta)\sigma(x,\theta)'$.
For estimation, we use the quasi log likelihood function
\begin{eqnarray*} 
{\mathbb H}_n(\theta) = 
-\frac{1}{2} \sum_{j=1}^n \left\{ \log \det S(X_{t_{j-1}}, \theta) 
+ h^{-1}  S^{-1}(X_{t_{j-1}},\theta) [ (\mathsf{\Delta}_j Y)^{\otimes 2}] \right\},
\end{eqnarray*}
where $\mathsf{\Delta}_jY=Y_{t_{j}}-Y_{t_{j-1}}$.
Then the quasi maximum likelihood estimator (QMLE) $\hat{\theta}^M_n$ is any estimator that satisfies 
\begin{eqnarray} \label{M-est}
\hat{\theta}^M_n &\in& \text{argmax }_{\theta \in \overline{\Theta}} {\mathbb H}_n(\theta).
\end{eqnarray}
The quasi Bayesian estimator (QBE) $\hat{\theta}^B_n$ 
with respect to the quadratic loss and 
a prior density $\pi : \Theta \rightarrow {{\mathbb R}}_+$ 
is given by 
\begin{eqnarray} \label{Bayes1}
\hat{\theta}^B_n = \left( \int_{\Theta} \exp ( {\mathbb H}_n(\theta) ) \pi(\theta) d\theta \right)^{-1}
\int_{\Theta} \theta \exp ( {\mathbb H}_n(\theta) ) \pi(\theta) d\theta. 
\end{eqnarray}
The prior density $\pi$ is assumed to be continuous and to satisfy 
$0 < \inf_{\theta \in \Theta} \pi(\theta) \leq \sup_{\theta \in \Theta} \pi(\theta) <\infty$.

For ${\bf n} =({ n}_1, \ldots, { n}_{\sf d})\in{\mathbb Z}_+^{\sf d}$ 
and ${\bf \nu} =({\bf \nu}_1, \ldots, {\bf \nu}_{{\sf p}})\in{\mathbb Z}_+^{\sf p}$, let 
$|{\bf n}| = { n}_1+ \ldots + { n}_{\sf d}$, 
$|{\bf \nu}| = {\bf \nu}_1+ \ldots + {\bf \nu}_{{\sf p}}$,
$\partial_x^{{\bf n}} = \partial_{x_1}^{{ n}_1} \cdots \partial_{x_{\sf d}}^{{ n}_{\sf d}}$, $\partial_{x_i} = \partial/\partial x_i$,
and $\partial_\theta^{\bf \nu} = \partial_{\theta_1}^{{\bf \nu}_1} \cdots \partial_{\theta_{{\sf p}}}^{{\bf \nu}_{{\sf p}}}$, 
$\partial_{\theta_i} = \partial/\partial \theta_i$.
We denote by $C_\uparrow^{k,l}({\mathbb R}^{\sf d} \times \Theta; {\mathbb R}^{{\sf m}})$ 
the set of all functions $f:{\mathbb R}^{\sf d} \times \Theta\to{\mathbb R}^{{\sf m}}$ satisfying the following conditions. 
(i) $f(x,\theta)$ is continuously differentiable with respect to $x$ up to order $k$
for all $\theta$.
(ii) For $|{\bf n}|=0,1, \ldots, k$, $\partial_x^{{\bf n}} f(x,\theta)$ is continuously differentiable
with respect to $\theta$ up to order $l$ for all $x$.
Moreover, for $|{\bf \nu}|=0,1, \ldots, l$ and $|{\bf n}|=0,1,\ldots, k$, 
$\partial_\theta^{\bf \nu} \partial_x^{\bf n} f(x,\theta)$ is of at most polynomial growth in $x$ uniformly in $\theta$.
We denote by $\rightarrow^{d_s({\cal F})}$ the ${\cal F}$-stable convergence in distribution. 
Suppose that $\Theta$ has a Lipschitz boundary. 

The following Condition [H1] is Condition [H1$^\sharp$] of \cite{uchida2013quasi}. 
\begin{description}
\item[[H1\!\!\!]] {(i)} 
$\sup_{0 \leq t \leq T} \|b_t\|_p < \infty$ for all $p>1$.
\begin{description}
\item[(ii)] 
$\sigma \in C_\uparrow^{2,4}({\mathbb R}^{\sf d} \times \Theta; {\mathbb R}^{{\sf m}} \otimes {\mathbb R}^{\sf r})$
and $\inf_{x, \theta} \det S(x,\theta) > 0$. 
\item[(iii)] 
The process $X$ has a representation 
\begin{eqnarray*}
X_t &=& X_0 + \int_0^t \tilde{b}_s ds + \int_0^t a_{s} dw_s + \int_0^t \tilde{a}_{s} d\tilde{w}_s, 
\end{eqnarray*}
where 
$\tilde{b}$, $a$ and $\tilde{a}$ are 
progressively measurable processes taking values in ${\mathbb R}^{\sf d}$, 
${{\mathbb R}}^{\sf d} \otimes {{\mathbb R}}^{\sf r}$ and ${{\mathbb R}}^{\sf d} \otimes {{\mathbb R}}^{{\sf r}_1}$, respectively, 
and satisfy 
\begin{eqnarray*} 
 \|X_0\|_p+
\sup_{t\in[0,T]}(\|\tilde{b}_t\|_p +\|a_t\|_p+\|\tilde{a}_t\|_p)<\infty
\end{eqnarray*}
for every $p>1$. 
$\tilde{w}$ is an ${{\sf r}_1}$-dimensional Wiener process independent of $w$,
\end{description}
\end{description}
\vspace*{3mm}

Let 
\relsize{-0.5}
\begin{eqnarray*}
{\mathbb Y}(\theta) &=& -\frac{1}{2 T}\int_0^T \left\{ \log \left( \frac{\det S(X_t,\theta)}{\det S(X_t,\theta^*)} \right)
+\mbox{Tr} \left( S^{-1}(X_t,\theta) S(X_t,\theta^*) -I_d \right) \right\} dt.
\end{eqnarray*}
\relsize{0.5}
A key index $\chi_0$ is defined by 
\begin{eqnarray}\label{230806-1}
\chi_0 = \inf_{\theta \ne \theta^*} \frac{- {\mathbb Y}(\theta) }{|\theta -\theta^*|^2}. 
\end{eqnarray}
Non-degeneracy of $\chi_0$ plays an important role in the discussion. 
\begin{description}
\item[[H2\!\!\!]] For every $L>0$, there exists a constant $c_L $ such that
\begin{eqnarray*}
P\left[ \chi_0 \leq r^{-1} \right] &\leq& \frac{c_L}{r^{L}}
\end{eqnarray*}
for all $r>0$. 
\end{description}

Define the random field ${\mathbb Z}_n$ on ${\mathbb U}_n $ by 
\begin{eqnarray}\label{230808-1} 
{\mathbb Z}_n(u) = \exp \left\{ {\mathbb H}_n \left( \theta^* + \frac{1}{\sqrt{n}} u \right) - {\mathbb H}_n (\theta^*) \right\}
\end{eqnarray}
for $u \in {\mathbb U}_n $. 
Then following the proof of Theorem 3 of \cite{uchida2013quasi}, we see that 
Condition [H2] together with [H1] implies that 
for every $L>0$, 
\begin{eqnarray}\label{290121-1} 
P \left[ \sup_{u \in V_n(r)} {\mathbb Z}_n(u) \geq e^{-r^{2-\epsilon}} \right] \leq \frac{C_L}{r^L}
\qquad(r>0,\>n\in{\mathbb N})
\end{eqnarray}
for some constant $C_L$ and some $\epsilon\in(0,1)$. 
Thus Condition [A1] is fulfilled for $a_n = n^{-1/2}I_{{\sf p}\times {\sf p}}$ in the present situation.

Let 
\begin{eqnarray*}
\Gamma(\theta^*)[u,u] &=& \frac{1}{2 T} \int_0^T \mbox{Tr} \left( (\partial_\theta S) S^{-1}
(\partial_\theta S) S^{-1}(X_t,\theta^*)[u^{\otimes 2}] \right) dt, 
\end{eqnarray*}
Now we have 
\begin{theorem} \label{thm3}{\rm (Theorems 4 and 5 of \cite{uchida2013quasi})} 
Suppose that $[H1]$ and $[H2]$ are satisfied. Then, 
for $A=M$ and $B$, 
$\sqrt{n} ( \hat{\theta}_n^A -\theta^*) \rightarrow^{d_s({\cal F})} \Gamma(\theta^*)^{-1/2} \zeta$ and
$$
E \left[ {\sf f}(\sqrt{n} ( \hat{\theta}_n^A -\theta^*) ) \right] \rightarrow {\mathbb E} \left[ {\sf f}(\Gamma(\theta^*)^{-1/2} \zeta) \right]
$$
as $n \rightarrow \infty$ for all continuous functions ${\sf f}$ of at most polynomial growth, 
where $\zeta$ is a ${\sf p}$-dimensional standard Gaussian vector independent of ${\cal F}$. 
\end{theorem}
\vspace*{3mm}


Let 
\begin{eqnarray*}
\Delta_n[u] &=& \frac{1}{\sqrt{n}} \partial_\theta {\mathbb H}_n(\theta^*)[u] \\
&=& -\frac{1}{2 \sqrt{n}} \sum_{j=1}^n \Biggl\{ (\partial_\theta \log \det S(X_{t_{j-1}}, \theta^*) )[u] \\
&&\hspace{90pt}+ h^{-1} (\partial_\theta S^{-1})(X_{t_{j-1}},\theta^*) [u, (\mathsf{\Delta}_k Y)^{\otimes 2}] \Biggr\}, 
\\
\Gamma_n(\theta)[u,u] &=& - \frac{1}{n} \partial_\theta^2 {\mathbb H}_n(\theta)[u,u] \\
&=& \frac{1}{2 n} \sum_{j=1}^n \Biggl\{ (\partial_\theta^2 \log \det S(X_{t_{j-1}}, \theta) )[u^{\otimes 2}] \\&&\hspace{90pt}+ h^{-1} (\partial_\theta^2 S^{-1})(X_{t_{j-1}},\theta) )[u^{\otimes 2}, (\mathsf{\Delta}_k Y)^{\otimes 2}] \Biggr\}
\end{eqnarray*}
and 
\begin{eqnarray*}
r_n(u) &=& \int_0^1 (1-s) \left\{ \Gamma(\theta^*) - \Gamma_n(\theta^* + sn^{-1/2}u)\right\}[u,u] ds.
\end{eqnarray*}
Then, for $u\in{\mathbb R}^{\sf p}$ and large $n$, we have 
\begin{eqnarray*} 
{\mathbb Z}_n(u) &=& \exp \left( \Delta_n[u] -\frac{1}{2} \Gamma(\theta^*)[u,u] + r_n(u) \right). 
\end{eqnarray*}

\begin{lemma} \label{lem2}{\rm (Lemma 7 of \cite{uchida2013quasi})}
Assume $[H1]$. Then, for every $q>0$,
\begin{description}
\item[(i)] 
$
\sup_{n \in {\bf N}} E 
\left[ \left( \sqrt{n} \left| \Gamma_n (\theta^*) - \Gamma (\theta^*) \right| \right)^q \right] <\infty,
$
\item[(ii)] 
$
\sup_{n \in {\bf N}} E 
\left[ \left( \frac{1}{n} \sup_{\theta \in \Theta}  \left| \partial_\theta^3 {\mathbb H}_n (\theta) \right| \right)^q \right] <\infty.
$
\end{description}
\end{lemma}
\vspace*{3mm}
%

%

Then [A12](iv) is verified by using Lemma \ref{lem2}
for $\Gamma=\Gamma(\theta^*)$ 
under the Condition [A6].
It is not difficult to check [A10](iii) and (iv) with stability of the convergence 
if one follows the proof of Lemma 9 of \cite{uchida2013quasi}. 
The $L^p$ boundedness of $\{\Delta_n\}_n$ is obvious for all $p>1$.  
Almost sure positive definiteness of $\Gamma$ (i.e., [A10](ii)) and 
$L^p$ integrability of $\det\Gamma^{-1}$ follow from [H2]. 
Thus all the conditions in Condition [A10] are satisfied. 
$L^p$ integrability of $\Gamma$ is obvious, therefore $|\Gamma^{-1}|$ is $L^p$ integrable, 
which implies Conditions [A12](ii) and (iii). 
Thus all the conditions in [A12] are satisfied.
Condition [A7'] holds obviously under [A12]. 
Condtions [A2-6], [A11] and [A13] are fulfilled easily for some $\Psi_n\in{\rm GL}(|\jz|)$. (See Example \ref{ex_bridge} in Section \ref{scvs} for instance.)
Consequently, the results in Sections \ref{spqle}-\ref{smconv} 
about the penalized estimators for (\ref{290121-5})
are valid. 

Condition [H2] can be easily verified if we apply the analytic criterion or the geometric criterion 
of \cite{uchida2013quasi}. \\


\section{Simulation study} \label{ssim}
In this section, we report the resutl of the simulation study to check the performance of the variable selection based on our penalized method.
The model is a volatility regression model in section \ref{sapp}.
Let ${\sf p}=d$, $a_n = n^{-1/2}I_{{\sf p}\times {\sf p}}$, $q < 1$, $\sigma(\theta,x) = \exp \bigl(\sum_{k=1}^{\sf p} \theta_k\sin(x_s^k)\bigr)$ and 
\begin{align}
	X^k_T = \int_0^t\frac{\sin(2k\pi s)}{(1+(X^k_s)^2)}dw_s^k \quad k = (1,\ldots,d), \notag
\end{align}
where $w^1,\ldots,w^d$ are independent standard Brownian motions.\\

\begin{table}[t]
\begin{center}
\caption{Simulation results for the volatility regression model.}\label{tab1}
\scalebox{0.847}{
\begin{tabular}{ccccccc}\hline
	&True&&1000&2000&3000&10000\\ \hline
	&&QMLE &0.023(0.260)&-0.006(0.181)&0.007(0.152)&0.001(0.083)\\
	$\hat{\theta}_1$&0&p-QL &-0.005(0.081)&0.002(0.051)&0.000(0.022)&0.000(0)\\
	&&prob&0.981&0.989&0.997&1\\
	&&QMLE &0.989(0.271)&1.002(0.180)&0.999(0.148)&1.000(0.085)\\
	$\hat{\theta}_2$&1&p-QL &0.792(0.396)&0.893(0.244)&0.937(0.160)&0.972(0.077) \\
	&&prob&0.869&0.971&0.996&1\\
	&&QMLE&0.003(0.263)&-0.004(0.181)&-0.009(0.146)&0.007(0.084) \\
	$\hat{\theta}_3$&0&p-QL&0.003(0.088)&0.000(0.030)&-0.002(0.042)&0.000(0) \\
	&&prob&0.982&0.992&0.995&1\\
	&&QMLE&0.986(0.263)&1.007(0.176)&1.003(0.146)&0.998(0.083) \\
	$\hat{\theta}_4$&1&p-QL&0.808(0.392)&0.898(0.241)&0.932(0.160)&0.968(0.073) \\
	&&prob&0.88&0.977&0.993&1\\
	&&QMLE&1.997(0.258)&2.000(0.184)&2.006(0.148)&1.996(0.084) \\
	$\hat{\theta}_5$&2&p-QL&1.912(0.328)&1.940(0.231)&1.965(0.139)&1.980(0.076) \\
	&&prob&0.999&0.999&1&1\\
	&&QMLE&-0.010(0.272)&0.001(0.185)&-0.004(0.155)&-0.001(0.083) \\
	$\hat{\theta}_6$&0&p-QL&0.001(0.123)&-0.003(0.059)&0.000(0.028)&0.000(0) \\
	&&prob&0.968&0.989&0.995&1\\
	&&QMLE&0.997(0.262)&0.996(0.179)&0.998(0.153)&0.998(0.081) \\
	$\hat{\theta}_7$&1&p-QL&0.789(0.390)&0.892(0.246)&0.927(0.172)&0.971(0.078) \\
	&&prob&0.867&0.967&0.992&1\\
	&&QMLE&1.000(0.267)&0.991(0.187)&1.002(0.150)&1.000(0.818) \\
	$\hat{\theta}_8$&1&p-QL&0.811(0.399)&0.883(0.248)&0.936(0.162)&0.972(0.076) \\
	&&prob&0.881&0.968&0.995&1\\
	&&QMLE&1.006(0.269)&0.997(0.182)&1.002(0.152)&0.999(0.083) \\
	$\hat{\theta}_9$&1&p-QL&0.788(0.394)&0.891(0.241)&0.937(0.167)&0.973(0.071) \\
	&&prob&0.871&0.971&0.994&1\\
	&&QMLE&0.022(0.264)&-0.008(0.181)&0.007(0.145)&0.000(0.081) \\
	$\hat{\theta}_{10}$&0&p-QL&0.000(0.089)&-0.004(0.048)&0.000(0.027)&0.000(0.007) \\
	&&prob&0.975&0.986&0.993&0.998\\ \hline
	&&Under model&0.91&0.958&0.98&0.998\\
	Total&&Over model&0.606&0.9&0.979&1\\
	&& True model&0.591&0.878&0.962&0.998\\ \hline
\end{tabular}
}
\end{center}
\end{table}

Obviously, Condition [H1] is fulfilled.
Following the section \ref{sapp}, we define ${\mathbb H}_n$ by
\begin{eqnarray} 
{\mathbb H}_n(\theta) = 
-\frac{1}{2} \sum_{j=1}^n \left\{ \log \det S(X_{t_{j-1}}, \theta) 
+ h^{-1}  S^{-1}(X_{t_{j-1}},\theta) [ (\mathsf{\Delta}_j Y)^{\otimes 2}] \right\}. \notag
\end{eqnarray}
Similarly to the proof of Theorems 5 of \cite{uchida2013quasi}), we have [A1].
Moreover we have [A10] and [A13] as discussed in section \ref{sapp}.
Define $p$ by $p(x) = |x|^q$ and $\xi_n^j$ by $\xi_n^j = n^{q'/2}$.
By definition, we have [A2-6], [A11] and [A14].
We set $q=0.3$, $q'=2/3$, ${\sf p} = d = 10$ and $T=1$.
The true value $\theta^*$ of an unknown parameter $\theta$ is $\theta^* = (0,1,0,1,2,0,1,1,1,0)'$.
Four cases of n are considered: $n = 1000, 2000, 3000$ and $10000$.
We used the local quadratic approximation in \cite{fan2001variable} for the optimization of the penalized quasi likelihood function.

Table \ref{tab1} compares the averages and standard deviations (parentheses) of quasi maximum likelihood estimator (QMLE) and penalized estimator (p-QL) over 1000 iterations for each cases.
Table \ref{tab1} also shows the probability that correct model is selected is selected:
\begin{align}
	P(\hat{\theta}_{n.j}=0)  \notag
\end{align}
for $j \in \jz$ and
\begin{align}
	P(\hat{\theta}_{n.j}\neq0)  \notag
\end{align}
for $j\in\jo$.
Under model is the probability that the estimator selects an under model:
\begin{align}
	P\Bigl(\{j;\hat{\theta}_{n,j}=0\} \supset\jz \Bigr) \notag
\end{align}
and Over model is the probability that the estimator selects an over model:
\begin{align}
	P\Bigl(\{j;\hat{\theta}_{n,j}=0\} \subset\jz \Bigr). \notag
\end{align}
True model is the probability that the true model is selected:
\begin{align}
	P\Bigl(\{j;\hat{\theta}_{n,j}=0\} =\jz \Bigr). \notag
\end{align}

From Table \ref{tab1}, it can seen that when sample size $n$ is large , penalized method performs variable selection very well. Moreover, the bias of non-zero parameters decrease as the sample size increases. 

\bibliographystyle{myspbasic}
\bibliography{bibtex}

\end{document}